\documentclass[11pt]{amsart}
\usepackage{amssymb,amscd,  latexsym, graphicx, mathrsfs, enumerate}

\newcommand{\nc}{\newcommand}
\numberwithin{equation}{section}
\newtheorem{theorem}{Theorem}[section]
\newtheorem{prop}[theorem]{Proposition}
\newtheorem{importnota}[theorem]{Important Notation}
\newtheorem{prblm}[theorem]{Problem}
\newtheorem{notation}[theorem]{Notation}
\newtheorem{caution}[theorem]{Caution}
\newtheorem{remark}[theorem]{Remark}
\newtheorem{lemma}[theorem]{Lemma}
\newtheorem{construction}[theorem]{Construction}
\newtheorem{corollary}[theorem]{Corollary}
\newtheorem{example}[theorem]{Example}
\newtheorem{conclusion}[theorem]{Conclusion}
\newtheorem{triviality}[theorem]{Triviality}
\newtheorem{proto}[theorem]{Prototype Quasifibration}
\newtheorem{cauex}[theorem]{Cautionary Example}
\newtheorem{propositiondef}[theorem]{Proposition-Definition}
\newtheorem{subth}{Nuisance}[theorem]
\newtheorem{ssubth}{ }[subth]
\newtheorem{conjecture}[theorem]{Conjecture}
\newtheorem{sidest}[theorem]{Side Story}
\newtheorem{miniexample}[theorem]{Example}
\theoremstyle{definition}
\newtheorem{defin}[theorem]{Definition}

\nc\tri[1]{\begin{triviality}}
\nc\side[1]{\begin{sidest}}
\nc\conj[1]{\begin{conjecture}}
\nc\prodef[1]{\begin{propositiondef}}
\nc\prt[1]{\begin{proto}}
\nc\lem[1]{\begin{lemma}}
\nc\sblm[1]{\begin{sublemma}}
\nc\pro[1]{\begin{prop}}
\nc\thm[1]{\begin{theorem}}
\nc\cor[1]{\begin{corollary}}
\nc\dfn[1]{\begin{defin}}
\nc\sthm[1]{\begin{subth}}
\nc\exm[1]{\begin{example}}
\nc\miniexm[1]{\begin{miniexample}}
\nc\plm[1]{\begin{prblm}}
\nc\rmk[1]{\begin{remark}}
\nc\subrmk[1]{\begin{subremark}}
\nc\ntn[1]{\begin{notation}}
\nc\cau[1]{\begin{caution}}
\nc\imn[1]{\begin{importnota}}
\nc\cax[1]{\begin{cauex}}
\nc\con[1]{\begin{construction}}
\nc\ssthm[1]{\begin{ssubth}}
\nc\cnc[1]{\begin{conclusion}}
\nc\elem{\end{lemma}}
\nc\esblm{\end{sublemma}}
\nc\eside{\end{sidest}}
\nc\econj{\end{conjecture}}
\nc\eprodef{\end{propositiondef}}
\nc\eprt{\end{proto}}
\nc\ethm{\end{theorem}}
\nc\ecor{\end{corollary}}
\nc\edfn{\end{defin}}
\nc\esthm{\end{subth}}
\nc\epro{\end{prop}}
\nc\etri{\end{triviality}}
\nc\eexm{\end{example}}
\nc\eminiexm{\end{miniexample}}
\nc\ermk{\end{remark}}
\nc\subermk{\end{subremark}}
\nc\eplm{\end{prblm}}
\nc\ecau{\end{caution}}
\nc\ecax{\end{cauex}}
\nc\eimn{\end{importnota}}
\nc\entn{\end{notation}}
\nc\econ{\end{construction}}
\nc\ecnc{\end{conclusion}}
\nc\essthm{\end{ssubth}}

\newcommand{\Q}{\mathbb{Q}}

\newcommand{\Z}{\mathbb{Z}}

\newcommand{\F}{\mathbb{F}}

\renewcommand{\Bbb}{\mathbb}

\title[Central limit theorem for Artin L-functions]
{Central limit theorem for Artin L-functions}

\author[Peter Cho]{Peter J. Cho}
\address{Department of
Mathematics, State University of New York at Buffalo, NY 14260, USA}
\email{jcho23@buffalo.edu}
\author[Henry Kim]{Henry H. Kim$^{\star}$}
\thanks{$^{\star}$ partially supported by an NSERC grant.}
\address{Department of
Mathematics, University of Toronto, Toronto, ON M5S 2E4, CANADA \\
and Korea Institute for Advanced Study, Seoul, Korea}
\email{henrykim@math.toronto.edu}
\begin{document}
\begin{abstract} We show that the sum of the traces of Frobenius elements of Artin
$L$-functions in a family of $G$-fields satisfies the Gaussian distribution under certain counting conjectures.
We prove the counting conjectures for $S_4$ and $S_5$-fields. We also show central limit theorem for modular form $L$-functions with the trivial central character with respect to congruence subgroups as the level goes to infinity.
\end{abstract}
\maketitle

\section{Introduction}

Let $G$ be a finite group which is a transitive subgroup of a certain symmetry group $S_{d+1}$. A number field $K$ of degree $d+1$ is called a $G$-field if its Galois closure $\widehat{K}$ over $\Q$ is a $G$-Galois extension. For a $G$-field $K$, we attach the Artin L-function
$$
L(s,\rho,K)=\frac{\zeta_K(s)}{\zeta(s)}=\sum_{n=1}^\infty a_{\rho}(n)n^{-s},
$$
where $\rho$ is $d$-dimensional representation of $G$. Note that $-1\leq a_{\rho}(p)\leq d$.
If $G=S_{d+1}$, $\rho$ is the $d$-dimensional standard representation of $S_{d+1}$.
Let $L(X)^{r_2}$ be the set of $G$-fields $K$ with $ |d_K| <X$ and signature $(r_1,r_2)$.
In this paper, we restrict to the case $G=S_{d+1}$, and consider the sum $\sum_{p\leq x} a_{\rho}(p)$ in the family of $S_{d+1}$-fields, $L(X)^{r_2}$ and show that under the counting conjectures (\ref{estimate}) and (\ref{estimate1}),
it has the Gaussian distribution, namely, for a continuous real function $h$ on $\Bbb R$, if $\frac {\log X}{\log x}\to\infty$ as $x\to\infty$,

\begin{equation}\label{main-id}
\frac 1{\#L(X)^{r_2}} \sum_{L(s,\rho,K)\in L(X)^{r_2}} h\left(\frac {\sum_{p\leq x} a_{\rho}(p)}{\sqrt{\pi(x)}}\right)
\longrightarrow \frac 1{\sqrt{2\pi}}\int_{-\infty}^\infty h(t)e^{-\frac {t^2}2}\, dt.
\end{equation}
When $L(X)^{r_2}$ is replaced by $\mathcal F_k$, the family of all normalized holomorphic Hecke eigen cusp forms of weight $k$ with respect to $SL_2(\Bbb Z)$,
Nagoshi \cite{N} showed that if $\frac {\log k}{\log x}\to\infty$ as $x\to\infty$,

\begin{equation}\label{N}
\frac 1{\#\mathcal F_k} \sum_{f\in \mathcal F_k} h\left(\frac {\sum_{p\leq x} a_{f}(p)}{\sqrt{\pi(x)}}\right)
\longrightarrow \frac 1{\sqrt{2\pi}}\int_{-\infty}^\infty h(t)e^{-\frac {t^2}2}\, dt,
\end{equation}
and called it central limit theorem.

Here we note that we do not need the Artin conjecture nor the strong Artin conjecture in the proof of (\ref{main-id}).
The estimates (\ref{estimate}) and (\ref{estimate1}) are proved by Taniguchi and Thorne \cite{TT} for $S_3$ fields.
For $G=S_4, S_5$, the estimate (\ref{estimate}) was proved in \cite{BBP} and \cite{ST}, resp.
We prove (\ref{estimate1}) in Sections \ref{S_5} and \ref{S_4}. Hence (\ref{main-id}) is unconditional for $S_3, S_4$ and $S_5$-fields.
These estimates will be used in computing the $n$-level densities of Artin $L$-functions \cite{CK}, \cite{CK1}.

We also study the distribution of the prime sums
$\sum_{p\leq x} a_{\rho}(p)^r$ for a positive integer $r$. The effective Chebotarev density theorem implies the analogue of Sato-Tate distribution. Namely, $\frac 1{\pi(x)}\sum_{p\leq x} a_{\rho}(p)^r\to n_r$ as $x\to\infty$, where $n_r$ is the multiplicity of the trivial representation in $\rho^r$.

We also study distribution of $\sum_{p\leq x} a_{f}(p)$ for $f\in S_k(N)$, the set of normalized Hecke eigen cusp forms of weight $k$ with respect to $\Gamma_0(N)$ with the trivial central character. We prove the Gaussian distribution as $N\to\infty$.

\section{Counting $S_{d+1}$-fields with local conditions}

Let $\mathcal S = (LC_{p})$ be a finite set of local
conditions. $LC_p=\mathcal S_{p,C}$ means that $p$ is unramified in $K$ and the conjugacy class of $\text{Frob}_p$ is $C$. Let $|\mathcal S_{p,C}|=\frac {|C|}{|G|(1+f(p))}$ for some positive function $f(p)$ which satisfies
$f(p)=O(\frac 1p)$. There are also several splitting types of ramified primes, which are denoted by $r_1,r_2,\cdots,r_w$. If $LC_p=\mathcal S_{p,r_i}$, then it means that $p$ is ramified and its splitting type is $r_i$. Assume that we can choose explicit positive functions $c_1(p),c_2(p), \cdots, c_w(p)$ with $\sum_{i=1}^w c_i(p)=f(p)$. Define $|\mathcal{S}_{p,r_i}|=\frac{c_i(p)}{1+f(p)}$ and $|\mathcal S|=\prod_p |LC_p|$.

Let $L(X;\mathcal S)^{r_2}$ be the set of $S_{d+1}$-fields $K$ with $ |d_K| <X$ and the local condition
$\mathcal S$.

Then assume that
\begin{eqnarray} \label{estimate}
|L(X)^{r_2}| &=& A(r_2) X +O(X^{\delta}),\\
|L(X;\mathcal S)^{r_2}| &=& |\mathcal S| A(r_2) X + O\left(\left(\prod_{p \in S} p \right)^{\gamma} X^{\delta} \right),\label{estimate1}
\end{eqnarray}
for some positive constant $\delta$ and $\gamma$, and the implied constant is uniformly bounded for $p$ and local conditions at $p$.

This assumption is satisfied when $G=S_3, S_4$ and $S_5$. When $G=S_3$, Taniguchi and Thorne \cite{TT} obtained more precise results:
Let $L(X)^\pm$ be the set of cubic fields $K$ with $ \pm d_K <X$. Then
\begin{equation*}
|L(X)^\pm|=\frac{C^{\pm}}{12 \zeta(3)}X + \frac{4 K^{\pm}}{5 \Gamma(2/3)^3}X^{5/6}+O(X^{7/9+\epsilon}),
\end{equation*}
where $C^{-}=3$, $C^+=1$, $K^-=\sqrt{3}$, and $K^+=1$. Here, we count only one cubic field from three conjugate fields. Let $TS_p, PS_p,$ and $IN_p$ be the local conditions of $p$ which means that $p$ is totally split, partially split and inert respectively. Let $S=\{LC_{p_i}| i=1,2,\cdots,u \}$ be a set of local conditions at $p_i$. Then
\begin{eqnarray*}
|LC_p|=\left\{ \begin{array}{cc}
\frac{1/6}{1+1/p+1/p^2} & \mbox{if $LC_p=TS_p$},\\
\frac{3/6}{1+1/p+1/p^2} & \mbox{if $LC_p=PS_p$},\\
\frac{2/6}{1+1/p+1/p^2} & \mbox{if $LC_p=IN_p$},\\
\frac{1/p}{1+1/p+1/p^2} & \mbox{if $p$ is partially ramified},\\
\frac{1/p^2}{1+1/p+1/p^2} & \mbox{if $p$ is totally ramified}.
\end{array}\right.
\end{eqnarray*}
and $|L(X; \mathcal S)^{\pm}| = |\mathcal S| A^\pm X + O(E_{\mathcal S}(X))$, where $A^\pm=C_1\frac{C^{\pm}}{12 \zeta(3)}$, and
\begin{eqnarray*}
E_{\mathcal S}(X)=\left\{ \begin{array}{lcc}
X^{5/6} & \mbox{ if $(\prod_{ p \in \mathcal S} p^{\frac{8}{9}e_p}) < X^{\frac{1}{16}}$, }\\
( \prod_{p \in \mathcal S}p^{\frac{8}{9}e_p})X^{\frac 79+\epsilon} & \mbox{ if $(\prod_{ p \in \mathcal S} p^{\frac{8}{9}e_p}) \geq X^{\frac{1}{16}}$, }\
\end{array}\right.
\end{eqnarray*}
where $e_p=1$ if $p$ is unramified and otherwise, $e_p=2$.

For $G=S_4, S_5$, the estimate (\ref{estimate}) was proved in \cite{BBP} and \cite{ST}.
We prove (\ref{estimate1}) in Sections \ref{S_5} and \ref{S_4}.

We identify $L(X)^{r_2}$ with the set of Artin $L$-functions $L(s,\rho,K)$ where $K$ is a $G$-field with $|d_K| <X$ with signature $(r_1,r_2)$. Here we count only one $G$-field for each $d+1$ conjugate fields.
Throughout the article, we
implicitly assume that the size of $L(X)^{r_2}$ is same with the number of number fields which
satisfy the requirements of $L(X)^{r_2}$. This claim deals with arithmetic equivalence of number
fields. Two number fields $K_1$ and $K_2$ are arithmetically equivalent if $\zeta_{K_1}(s) = \zeta_{K_2}(s)$. We say
that a number field K is arithmetically solitary if $\zeta_K(s) = \zeta_F(s)$ implies that $K$ and $F$ are
conjugate. It is known that $S_d$-fields and $A_d$ fields are arithmetically solitary. See Chapter
II in \cite{Kl}.

For simplicity, we denote $L(s,\rho,K)\in L(X)^{r_2}$ by
$\rho\in L(X)$.

\section{Central limit theorem of Artin L-functions}

Consider, for a continuous real function $h$ on $\Bbb R$,

\begin{equation}\label{central}
\frac 1{\#L(X)} \sum_{\rho\in L(X)} h\left(\frac {\sum_{p\leq x} a_{\rho}(p)}{\sqrt{\pi(x)}}\right).
\end{equation}
We assume that $x$ grows more slowly than $X$; namely, $\frac {\log X}{\log x}\longrightarrow \infty$ as $x\to\infty$. So for an arbitrary positive real number $a$, we have $X>x^a$.

By Theorem 25.8 and Theorem 30.2 (the method of moments) in \cite{B}, it is enough to consider $h(x)=x^r$.
Consider
\begin{equation}\label{central-limit}
\sum_{\rho\in L(X)} \left(\frac {\sum_{p\leq x} a_\rho(p)}{\sqrt{\pi(x)}}\right)^r.
\end{equation}

By multinomial formula,
$$\left(\sum_{p\leq x} a_{\rho}(p)\right)^r=\sum_{u=1}^r {\sum}^{(1)}_{(r_1,...,r_u)} \frac {r!}{r_1!\cdots r_u!} \frac 1{u!} {\sum}^{(2)}_{(p_1,...,p_u)} a_{\rho}(p_1)^{r_1}\cdots a_{\rho}(p_u)^{r_u},
$$
where $\sum_{(r_1,...,r_u)}^{(1)}$ means the sum over the $u$-tuples $(r_1,...,r_u)$ of positive integers such that $r_1+\cdots+r_u=r$, and $\sum_{(p_1,...,p_u)}^{(2)}$ means the sum over the $u$-tuples $(p_1,...,p_u)$ of distinct primes such that $p_i\leq x$ for each $i$.
Then
$$(\ref{central-limit})=\pi(x)^{-\frac r2} \sum_{u=1}^r \frac 1{u!} {\sum}^{(1)}_{(r_1,...,r_u)} \frac {r!}{r_1!\cdots r_u!}  {\sum}^{(2)}_{(p_1,...,p_u)} \left(\sum_{\rho\in L(X)} a_{\rho}(p_1)^{r_1}\cdots a_{\rho}(p_u)^{r_u}\right).
$$

Now we claim that except when $r$ is even, $u=\frac r2$, and $r_1=\cdots=r_u=2$, it gives rise to the error term.

Now suppose $r_i\geq 2$ for all $i$, and $r_j>2$ for some $j$. Then since $r_1+\cdots+r_u=r$, $u\leq\frac {r-1}2.$ Hence by the trivial estimate,
such term is majorized by
$$\pi(x)^{-\frac r2} \sum_{u=1}^r \frac 1{u!} {\sum}_{(r_1,...,r_u)}^{(1)}\frac {r!}{r_1!\cdots r_u!} d^{r_1+\cdots+r_u} |L(X)| \pi(x)^u
\ll_r \pi(x)^{-\frac 12} |L(X)| \sum_{u=1}^r \frac 1{u!}  d^r \ll_{r,d} X \pi(x)^{-\frac 12}.
$$
This gives rise to the error term.

Suppose $r_i\leq 2$ for all $i$.
Suppose $r_i=1$ for some $i$. We may assume that $r_1=1$.

Let $N$ be the number of conjugacy classes of $G$, and partition the sum $\sum_{\rho\in L(X)}$ into $(N+w)^{u}$ sums, namely, given
$(\mathcal S_1,...,\mathcal S_{u})$, where $\mathcal S_i$ is either $\mathcal S_{p_i,C}$ or $\mathcal S_{p_i,r_j}$,
we consider the set of $\rho\in L(X)$ with the local conditions $\mathcal S_i$ for each $i$. Note that in each such partition, $a_{\rho}(p_1)^{r_1}\cdots a_{\rho}(p_u)^{r_u}$ remains a constant.

Suppose $p_1$ is unramified, and fix the splitting types of $p_2,\cdots,p_u$, and let $\text{Frob}_{p_1}$ runs through the conjugacy classes of $G$. Then by (\ref{estimate1}), the sum of such $N$ partitions is
$$\sum_C \left(\frac{|C|a_\rho(p_1)}{|G|(1+f(p_1))} A(\mathcal S_2,...,\mathcal S_{u})X + O((p_1\cdots p_u)^\gamma X^\delta) \right),
$$
for a constant $A(\mathcal S_2,...,\mathcal S_u)$.
Let $\chi_\rho$ be the character of $\rho$. Then $a_{\rho}(p)=\chi_{\rho}(g)$, where $g=\text{Frob}_p$. By orthogonality of characters,
$\sum_C |C| a_{\rho}(p_1)=\sum_{g\in G} \chi_\rho(g)=0$. Hence the above sum is
$O((p_1\cdots p_u)^\gamma X^\delta)$, and it is majorized by $\pi(x)^{-\frac r2+u}x^{\gamma u} X^{\delta}.$

Hence we can assume that $r_i\leq 2$ for each $i$, and $p_j$ is ramified when $r_j=1$. Suppose $r_1+\cdots+r_{v}+r_{v+1}+\cdots+r_u=r$, $r_1=\cdots=r_v=1$ and $r_{v+1}=\cdots=r_u=2$. Then $u-v\leq \frac {r-1}2$, and $p_1,...,p_v$ are ramified.
The partition of fixed splitting types of $p_{v+1},...,p_u$ is majorized by
$$
\prod_{i=1}^v \frac {f(p_i)}{1+f(p_i)} B(\mathcal S_{v+1},...,\mathcal S_u)X + O((p_1\cdots p_u)^\gamma X^\delta),
$$
for some constant $B(\mathcal S_{v+1},...,\mathcal S_u)$.
Since $\frac {f(p)}{1+f(p)}\ll \frac 1p$, it contributes to
$$\pi(x)^{u-v-\frac r2}(\log\log x)^v X+ \pi(x)^{-\frac r2+u}x^{\gamma u} X^{\delta}\ll  X (\log\log x)^v \pi(x)^{-\frac 12}+\pi(x)^{-\frac r2+u}x^{\gamma u} X^{\delta}.
$$

Now let $r$ be even, $u=\frac r2$, and $r_1=\cdots=r_u=2$. If one of $p_1,p_2, \cdots , p_u$ is ramified, their contribution is
majorized by $X \pi(x)^{-1}\log\log x.$
Now we assume that all primes are unramified. Then the corresponding term is
\begin{equation}\label{main}
\pi(x)^{-\frac r2} \frac 1{u!} \frac {r!}{2^u} \sum_{(p_1,...,p_u)}^{(2)} \left(\sum_{L(s,\rho)\in L(X)} a_{\rho}(p_1)^2\cdots a_{\rho}(p_u)^2\right).
\end{equation}

Let $N$ be the number of conjugacy classes of $G$, and partition the sum $\sum_{\rho\in L(X)}$ into $N^u$ sums
where $(C_1,...,C_u)$ is the set of $\rho\in L(X)$ such that $\text{Frob}_{p_i}\in C_i$ for each $i$.
Then,
\begin{eqnarray*}
&& \sum_{\rho\in L(X)} a_{\rho}(p_1)^2\cdots a_{\rho}(p_u)^2=\sum_{(C_1,...,C_u)} \chi_{\rho}(p_1)^2\cdots \chi_{\rho}(p_u)^2
\left(\sum_{\rho\in L(X)\atop \text{Frob}_{p_i}\in C_i} 1\right) \\
&& =\sum_{(C_1,...,C_u)} \chi_{\rho}(p_1)^2\cdots \chi_{\rho}(p_u)^2 \left(\prod_{i=1}^u \frac {|C_i|}{|G|(1+f(p_i))} |L(X)|+O((p_1\cdots p_u)^\gamma X^{\delta})\right).
\end{eqnarray*}
Now
$$
\sum_{(C_1,...,C_u)} \chi_{\rho}(p_1)^2\cdots \chi_{\rho}(p_u)^2 \prod_{i=1}^u \frac {|C_i|}{|G|(1+f(p_i))}=\prod_{i=1}^u \left(\sum_{C_i} \frac {\chi_{\rho}(p_i)^2|C_i|}{|G|(1+f(p_i))}\right).
$$
Here $\chi_{\rho}(p)^2=\chi_{\rho^2}(p)=\chi_{Sym^2\rho}(p)+\chi_{\wedge^2\rho}(p)$. We observed in \cite{CK} that since $\rho$ is an irreducible real self-dual representation,
$Sym^2\rho$ contains the trivial representation and
$\wedge^2\rho$ does not contain the trivial representation (\cite{JL}, page 274). Hence $\chi_{\rho}(p)^2=1+\sum_{j=1}^l \eta_j(p)$, where $\eta_j$'s are non-trivial irreducible characters of $G$. By the orthogonality of characters, for each $j$, $\sum_C |C|\eta_j(p)=\sum_{g\in G} \eta_i(g)=0$. Hence $\sum_{C} \chi_{\rho}(p)^2|C|=|G|.$
Therefore,
\begin{eqnarray*}
&& {\sum}_{(p_1,...,p_u)}^{(2)} \left(\sum_{\rho\in L(X)} a_{\rho}(p_1)^2\cdots a_{\rho}(p_u)^2\right) \\
&=& \pi(x)^u |L(X)| + O(\pi(x)^{u-1} |L(X)|\log\log x) + O(\pi(x)^u x^{\gamma u} X^{\delta}).
\end{eqnarray*}

Note
$$\frac{1}{\sqrt{2\pi}}\int_{-\infty}^{\infty} t^r e^{-\frac{t^2}{2}}dt = \begin{cases} \frac{r!}{(r/2)! 2^{r/2}}, & \text{if $r$ is even,}\\
0, & \text{if $r$ is odd}\end{cases}.
$$
Hence we have proved
\begin{theorem}\label{Artin}
Suppose $\frac {\log X}{\log x}\longrightarrow \infty$ as $x\to\infty$. Then
$$
\frac 1{|L(X)|} \sum_{\rho\in L(X)} \left(\frac {\sum_{p\leq x} a_{\rho}(p)}{\sqrt{\pi(x)}}\right)^r
=\frac{1}{\sqrt{2\pi}}\int_{-\infty}^{\infty} t^r e^{-\frac{t^2}{2}}dt + O\left(\frac {(\log\log x)^r}{\pi(x)^{\frac 12}}\right).
$$
\end{theorem}

This proves (\ref{main-id}).

\section{Central Limit Theorem for Hecke eigenforms; Level aspect}

In this section, in analogy to (\ref{N}),
we consider central limit theorem for modular form $L$-functions with the trivial central character with respect to congruence subgroups as the level goes to infinity. We follow \cite{N} closely. For $k\geq 2$, let $S_k(N)$ be the set of normalized Hecke eigen cusp forms of weight $k$ with respect to $\Gamma_0(N)$ with the trivial central character.
Let $f(z)=\sum_{n=1}^\infty a_f(n)n^{\frac {k-1}2} e^{2\pi inz}$; $a_f(mn)=a_f(m)a_f(n)$, if $(m,n)=1$; $a_f(1)=1$; $a_f(p^j)=a_f(p)a_f(p^{j-1})-a_f(p^{j-2})$.

We show
\begin{theorem} \label{Hecke}
For a continuous real function $h$ on $\Bbb R$, (assume that $\frac {\log N}{\log x}\longrightarrow \infty$ as $x\to\infty$.)

$$
\frac 1{\#S_k(N)} \sum_{f\in S_k(N)} h\left(\frac {\sum_{p\leq x} a_f(p)}{\sqrt{\pi(x)}}\right)\longrightarrow \frac 1{\sqrt{2\pi}}\int_{-\infty}^\infty h(t) e^{-\frac {t^2}2}\, dt\quad \text{as $x\to\infty$}.
$$
\end{theorem}

We have, from \cite{Se},

\begin{lemma} Suppose $k\geq 2$. Let $S_k(N,\chi)$ be the set of normalized Hecke eigen cusp forms of weight $k$ with respect to $\Gamma_0(N)$ with a character $\chi$ (mod $N$). Then
$$\sum_{f\in S_k(N,\chi)} a_f(n)=\frac {k-1}{12} \chi(\sqrt{n}) n^{-\frac 12} \psi(N)+O(n^{c}N^{\frac 12} d(N)),
$$
for some constant $c$, independent of $n, N$.
\end{lemma}

Here $\psi(N)=N \prod_{l | N} (1+\frac 1l)$,
and $d(N)$ is the number of positive divisors of $N$. Note that $\psi(N)=|SL_2(\Bbb Z): \Gamma_0(N)|$.
Here $\chi(x)=0$ if $x$ is not a positive integer prime to $N$. In particular, if $n$ is not a square,
$\sum_{f\in S_k(N,\chi)} a_f(n)=O(n^{c}N^{\frac 12} d(N)).$ Taking $n=1$ and $\chi=1$, we have
$$\#S_k(N)=\frac {k-1}{12} \psi(N)+O(N^{\frac 12}d(N)).
$$

We need to compute, for a positive integer $r$,
\begin{equation}\label{central-limit-h}
\sum_{f\in S_k(N)} \left(\frac {\sum_{p\leq x} a_f(p)}{\sqrt{\pi(x)}}\right)^r.
\end{equation}

By multinomial formula,
$$(\ref{central-limit-h})=\pi(x)^{-\frac r2} \sum_{u=1}^r \frac 1{u!} {\sum}_{(r_1,...,r_u)}^{(1)} \frac {r!}{r_1!\cdots r_u!} {\sum}_{(p_1,...,p_u)}^{(2)} \left(\sum_{f\in S_k(N)} a_f(p_1)^{r_1}\cdots a_f(p_u)^{r_u}\right).
$$

Now we claim that except when $r$ is even, $u=\frac r2$, and $r_1=\cdots=r_u=2$, it gives rise to the error term.

By \cite{N}, Lemma 2, we can show that
$a_f(p)^n=\sum_{j=0}^n h_n(j)a_f(p^j)$, where $h_n(j)=0$ if $n$ is odd and $j$ is even, or if $n$ is even and $j$ is odd.
For $u$-tuples $(r_1,...,r_u)$ and $(p_1,...,p_u)$, we define
\begin{eqnarray*}
A(r_1,...,r_u) &=& {\sum}_{(p_1,...,p_u)}^{(2)} B(r_1,...,r_u; p_1,...,p_u),\\
 B(r_1,...,r_u; p_1,...,p_u) &=& \sum_{f\in S_k(N)} a_f(p_1)^{r_1}\cdots a_f(p_u)^{r_u}.
\end{eqnarray*}

Then
$$
B(r_1,...,r_u; p_1,...,p_u)=\sum_{0\leq j_r\leq r_1,...,0\leq j_u\leq r_u} h_{r_1}(j_1)\cdots h_{r_u}(j_u) \sum_{f\in S_k(N)} a_f(p_1^{j_1}\cdots p_u^{j_u}).
$$

As in \cite{N}, if $r_l$ is odd for some $l$, $A(r_1,...,r_u)\ll N^{\frac 12}d(N) \pi(x)^u x^{cur}.$
	
Now let $r_1=\cdots=r_u=2$. Then $r$ is even, and $u=\frac r2$.

$$A(r_1,...,r_u)=\pi(x)^{\frac r2} \#S_k(N)+O(\pi(x)^{\frac r2-1}(\log\log x)^{\frac r2}\#S_k(N)).
$$

Now suppose that all $r_i$'s are even, and $r_i>2$ for some $i$. Then $u\leq \frac r2-1$.
Then
$$A(r_1,...,r_u)\ll \pi(x)^{-1}\#S_k(N).
$$

Hence, as in Theorem \ref{Artin}, we have
\begin{prop} Assume that $\frac {\log N}{\log x}\longrightarrow \infty$ as $x\to\infty$. Then
$$
\frac{1}{\#S_k(N)} \sum_{f\in S_k(N)} \left(\frac {\sum_{p\leq x} a_f(p)}{\sqrt{\pi(x)}}\right)^r
=\frac{1}{\sqrt{2\pi}}\int_{-\infty}^{\infty} t^r e^{-\frac{t^2}{2}}dt + O\left(\frac {(\log\log x)^{\frac r2}}{\pi(x)}\right).
$$
\end{prop}
This proves Theorem \ref{Hecke}

\section{Analogues of Sato-Tate distribution}

For a Hecke eigenform $f\in \mathcal F_k$, Sato-Tate conjecture says that for a continuous real function $h$ on $[-2,2]$,
$$
\frac 1{\pi(x)} \sum_{p\leq x} h(a_{f}(p))\longrightarrow \frac 1{2\pi}\int_{-2}^2 h(t) \sqrt{4-t^2}\, dt,\quad \text{as $x\to\infty$}.
$$
Let $a_f(p)=2\cos\theta_f(p)$ for $\theta_f(p)\in [0,\pi]$. Then
$\{\theta_f(p)\}$ is uniformly distributed with respect to the measure $\frac 2{\pi}\sin^2\theta \, d\theta$ on $[0,\pi]$. This is proved in \cite{BGHT}.

For a vertical Sato-Tate distribution, one can consider, for a fixed prime $p$,
\begin{equation} \label{CDF}
\sum_{f\in \mathcal F_k} a_{f}(p)^n.
\end{equation}

Conrey-Duke-Farmer \cite{CDF} proved, for a holomorphic form of weight $k$,
$$\sum_{f\in \mathcal F_k} a_f(p)^n=\frac k{6\pi} \left(1+\frac 1p\right)\int_0^\pi 2^n \cos^n\theta \frac {\sin^2\theta}{(1-\frac 1{p})^2+\frac 4p \sin^2\theta} \, d\theta+O(p^{\frac n2+\epsilon}).
$$

This implies that $\{\theta_f(p), f\in \mathcal F_k\}$ is uniformly distributed with respect to the measure
$$\frac 2{\pi} \left(1+\frac 1p\right)\frac {\sin^2\theta}{(1-\frac 1p)^2+\frac 4p \sin^2 \theta} \, d\theta.
$$

For Artin $L$-function analogue of Sato-Tate distribution, we consider, for $r\geq 1$,

\begin{equation}\label{Sato-Tate}
\frac 1{\pi(x)} \sum_{p\leq x} a_{\rho}(p)^r.
\end{equation}

In our case, note that $-1\leq a_{\rho}(p)\leq d$.
By effective Chebotarev density theorem (cf. \cite{Se1}, page 132), for $\log x\gg |G|(\log \left|d_{\widehat{K}}\right| )^2$,
\begin{equation*}\label{chebo}
\sum_{p\leq x \atop \text{Frob}_p\in C} 1=\frac {|C|}{|G|} \pi(x)+O\left(\pi(x^{\beta})\right)+O\left(x e^{-c |G|^{-\frac 12} (\log x)^{\frac 12}}\right),
\end{equation*}
where $\beta$ is an exceptional zero of $\zeta_{\widehat{K}}(s)$ such that $1-\beta\leq \frac 14\log d_{\widehat{K}}$, if it exists. Hence

$$\sum_{p\leq x} a_{\rho}(p)^r=\sum_C a_{\rho}(p)^r \left(\sum_{p\leq x\atop \text{Frob}_p\in C} 1\right)
=\sum_C a_{\rho}(p)^r \frac {|C|}{|G|} \pi(x)+ O(\pi(x^{\beta})+x e^{-c |G|^{-\frac 12} (\log x)^{\frac 12}}).
$$

Now $\sum_C |C| a_{\rho}(p)^r=\sum_{g\in G} \chi_{\rho}(g)^r$ and $\chi_{\rho}(g)^r=\chi_{\rho^r}(g)$.
Note that
$$\frac 1{|G|}\sum_{g\in G} \chi_{\rho^r}(g)=n_r,
$$
which is the multiplicity of the trivial representation in $\rho^r$. Hence
$$\sum_{p\leq x} a_{\rho}(p)^r=n_r \pi(x)+O(\pi(x^{\beta})+x e^{-c |G|^{-\frac 12} (\log x)^{\frac 12}}).
$$
Therefore,
$$\frac 1{\pi(x)} \sum_{p\leq x} a_{\rho}(p)^r\longrightarrow n_r, \quad \text{as $x\to\infty$}.
$$

For vertical Sato-Tate distribution, for a fixed prime $p$, consider
$$\frac 1{|L(X)|} \sum_{\rho\in L(X)} a_{\rho}(p)^r.
$$
Then by (\ref{estimate1}),
\begin{eqnarray*}
&& \sum_{\rho\in L(X)} a_{\rho}(p)^r=\sum_C a_{\rho}(p)^r \left(\sum_{\rho\in L(X)\atop \text{Frob}_p\in C} 1\right)+ a_{\rho}(p)^r
\left(\sum_{\rho\in L(X)\atop \text{$p$ is ramified}}  1 \right) \\
&&
=\frac {|L(X)|}{|G|(1+f(p))} \sum_C |C| a_{\rho}(p)^r +O(p^\gamma X^\delta)+O\left(\frac Xp\right)
=\frac {|L(X)| n_r}{1+f(p)}+O(p^\gamma X^\delta)+O\left(\frac Xp\right).
\end{eqnarray*}
So if $X>p^{\frac {1+\gamma}{1-\delta}}$,
$$\frac 1{|L(X)|} \sum_{\rho\in L(X)} a_{\rho}(p)^r=\frac {n_r}{1+f(p)}+O(p^{-1}).
$$

\section{Counting $S_5$ quintic fields with local conditions}\label{S_5}

Shankar and Tsimerman \cite{ST} recently counted $S_5$ quintic fields with a power saving error terms. For $i=0,1,2$, let $N_5^{(i)}(X)$ be the number of $S_5$ quintic fields of signature $(5-2i,i)$ with $|d_K| < X$. Then they showed
\begin{eqnarray*}
N_5^{(i)}(X)= D_i X+ O_{\epsilon}\left( X^{\frac{399}{400}+\epsilon }\right),
\end{eqnarray*}
where $D_i=d_i \prod_p (1+ p^{-2} -p^{-4} -p^{-5})$ and $d_0,d_1,d_2$ are $\frac{1}{240}, \frac{1}{24}$ and $\frac{1}{16}$, respectively.

We can count quintic fields with finitely many local conditions. Let $C$ be a conjugacy class of $S_5$ and $f(p)=p^{-1}+2p^{-2}+2p^{-3}+p^{-4}$. Let $\mathcal S=\{ LC_p \}$ be a finite set of local conditions. Define
 $|\mathcal S_{p,C}|=\frac {|C|}{|G|(1+f(p))}$, $|\mathcal S_{p,r_i} |=\frac {c_i(p)}{(1+f(p))}$, and $|\mathcal S|=\prod_p |LC_p|$, where $c_i(p)$'s are given explicitly at the end of this section.

\begin{theorem} Let $N_5^{(i)}(X,\mathcal S)$ be the number of $S_5$ quintic fields of signature $(5-2i,i)$ with $|d_K| < X$, and with the local condition $\mathcal S$. Then
\begin{eqnarray*}
N_5^{(i)}(X,\mathcal S)= |\mathcal S| D_i X+ O_{\epsilon}\left(\left(\prod_{p\in \mathcal S} p \right)^{2-\epsilon} X^{\frac{199}{200}+\epsilon }\right).
\end{eqnarray*}
\end{theorem}

We follow the notations in \cite{ST}.  Let $V_{\Z}$ be the space of $4$-tuples of $5 \times 5$ alternating matrices with integer coefficients. The group $G_\Z = GL_4(\Z) \times SL_5(\Z)$ acts on $V_\Z$ via
$$
(g_4,g_5) \cdot ( A, B, C, D)^t = g_4 ( g_5 A g_5^t, g_5 B g_5^t,g_5 C g_5^t,g_5 D g_5^t)^t.
$$
Here $g_4\cdot (A,B,C,D)^t$ means $(a_1 (A,B,C,D)^t, a_2 (A,B,C,D)^t, a_3 (A,B,C,D)^t, a_4 (A,B,C,D)^t)$, where $a_i$ is the $i$th row of $g_4$.

There is a canonical bijection between the set of $G_\Z$-equivalence classes of elements $(A,B,C,D) \in V_{\Z}$, and the set of isomorphism classes of pairs of $(R,R')$, where $R$ is a quintic ring and $R'$ is a sextic resolvent ring of $R$. (See \cite{B08}.)
Let $\mathcal V$ be an element of $V_{\Z}$.  Over the residue field $\F_p$, the element $\mathcal V$ determines a quintic
$\F_p$-algebra $R(\mathcal V)/(p)$. Let us define the splitting symbol $(\mathcal V,p)$ by
$$
(\mathcal V,p)=(f_1^{e_1}f_2^{e_2} \cdots ),
$$
whenever $R(\mathcal V)/(p) \cong \F_{p^{f_1}}[t_1]/(t_1^{e_1}) \oplus \F_{p^{f_2}}[t_2]/(t_2^{e_2}) \oplus \cdots$. Then there are 17 possible splitting types for $(\mathcal V,p)$; $(11111),$ $(1112),$ $(122),$ $(113),$ $(23),$ $(14),$ $(5),$ $(1^2111),$ $(1^212),$ $(1^23),$ $(1^21^21),$ $(2^21),$ $(1^311),$ $(1^32),$ $(1^31^2),$ $(1^41),$ and $(1^5)$. Let $\sigma$ be one of 17 splitting types. Then define $T_p(\sigma)$ to be the set of $\mathcal V\in V_{\Z}$ such that $(\mathcal V,p)=\sigma$ and $U_p(\sigma)$ to be the set of elements in $T_p(\sigma)$ corresponding to quintic rings that are maximal at $p$. The set $U_p(\sigma)$ is defined by congruence conditions on coefficients of $\mathcal V$ modulo $p^2$. Let $\mu(U_p(\sigma))$ be the $p$-adic density of $\mathcal S$ in $V_{\Z_p}$. They are computed in Lemma 4 in \cite{B08}. Let ${U}_p$ denote the union of the 17 $U_p(\sigma)$. Then Lemma 20 of \cite{B08} implies that
\begin{eqnarray*}
\mu({U_p})=(p-1)^8p^{12}(p+1)^4(p^2+1)^2(p^2+p+1)^2(p^4+p^3+p^2+p+1)(p^4+p^3+2p^2+2p+1)/p^{40}.
\end{eqnarray*}

Note that
$$d_i \zeta(2)^2\zeta(3)^2\zeta(4)^2\zeta(5) \prod_p \mu({U_p}) = d_i \prod_p (1+p^{-2}-p^{-4}-p^{-5}),
$$
which is the coefficient of the main term in counting quintic fields.  Here we need to interpret $\mu({U_p})$ in the following way:
$U_p$ can be considered as a subset of $(\Z/q^2 \Z)^{40}$, or the union of $k$ translates of $p^2 V_{\Z}$, where $k$ is the size of the set. Here $k$ is $\mu(U_p)q^{80}$.
Let ${W}_p$ be the complement of ${U_p}$ in $V_{\Z}$, then $\mu({W}_p)=1- \mu({U_p})$.  Then $W_p$ is the union of $\mu({W}_p)q^{80}$ translates of $q^2V_{\Z}$.

For $q$ square-free, let $W_q \subset V_{\Z}$ be the set of elements corresponding to quintic rings that are not maximal at each prime dividing $q$. Then $W_q$ is the union of $\prod_{p \mid q} \mu(W_p) \cdot q^{80}$ translates of $q^2 V_{\Z}$ by the Chinese Remainder Theorem.

Let $V_{\Z}^{ndeg}$ denote the set of elements in $V_{\Z}$ that correspond to orders in $S_5$-fields, and let $V_{\Z}^{deg}$ be the complement of $V_{\Z}^{ndeg}$.
A point in $V_{\Z}$ corresponds to a maximal order in an $S_5$ quintic fields precisely if it is in $ \cap_p U_p \cap V_{\Z}^{ndeg}$. For a $G_\Z$-invariant subset $S$ of $V_{\Z}$, let $N(S,X)$ denote the number of irreducible $G_\Z$-orbits in $S^{ndeg}:=S \cap V_{\Z}^{ndeg}$ having discriminant bounded by $X$.  For a set $S$ which is not $G_\Z$-invariant, we can define $N^*(S,X)$ which also counts the orbits of degenerate points in $S$.

Now we choose a finite set of primes $\{p_1,p_2, \cdots, p_n \}$. And choose a splitting type $\sigma_{p_k}$ for each $p_k$, $k=1,2,\cdots, n$. Define $U'_p$ to be $U_p$ if $p \neq p_k$, $k=1,2,\cdots, n$. If $p=p_k$ for some $k$, then $U'_p=U_p(\sigma_{p_k})$.   Let $W'_p$ be the complement of $U'_p$. Then $W'_q$ is the union of $\prod_{p \mid q} \mu(W'_p) \cdot q^{80}$ translates of $q^2 V_{\Z}$.

Let $N_5^{(i)}(X, \{ \sigma_{p_k} \}_{k=1,2,\cdots,n})$ be the number of $S_5$ quintic fields of signature $(5-2i,i)$ with $|d_K|<X$ and the splitting types of $p_k$'s are $\sigma_{p_k}$ for $k=1,2,\cdots,n$.
Then by inclusion-exclusion method,
$$N_5^{(i)}(X, \{ \sigma_{p_k} \}_{k=1,2,\cdots,n})=\sum_{q} \mu(q) N( W'_q \cap V_{\Z}^{(i)}, X).$$

We use two estimates of $N(W'_q \cap V_{\Z}^{(i)}, X)$ for small $q$'s and large $q$'s separately. Lemma 3 in \cite{ST} says that
$N(W_q, X)=O_\epsilon(X/q^{2-\epsilon})$ and it implies that $N(W'_q,X) \ll_\epsilon (p_1 p_2 \cdots p_n)^{2-\epsilon} \frac{X}{q^{2-\epsilon}}$ for $q$ square-free. If $L$ is a translate of $mV_{\Z}$, then we have
\begin{eqnarray}
N^*( \{ x \in L \cap V_{\Z}^{(i)} : a_{12} \neq 0\}, X)=c_i \frac{X}{m^{40}} + O\left(m^{-39}X^{\frac{39}{40}}\right),
\end{eqnarray}
where $c_i = d_i \zeta(2)^2\zeta(3)^2\zeta(4)^2\zeta(5)$. (See Equation 28 in \cite{B10}.)
Since $W'_q$ is an union of $\prod_{p \mid q} \mu(W'_p) \cdot q^{80}$ translates of $q^2 V_{\Z}$,
\begin{eqnarray*}
N^*( \{ x \in W'_q \cap V_{\Z}^{(i)} : a_{12} \neq 0\}, X)
=\prod_{p \mid q} \mu(W'_p) \cdot c_i X + O\left(q^2X^{\frac{39}{40}}\right)
 =\mu(W'_q)\cdot c_i X + O\left(q^2X^{\frac{39}{40}}\right).
\end{eqnarray*}

Then
\begin{eqnarray*}
&& N_5^{(i)}(X, \{ \sigma_{p_i} \}_{i=1,2,\cdots,n})=\sum_{q} \mu(q) N( W'_q \cap V_{\Z}^{(i)}, X)\\
&&= \sum_{q \leq Q} \left( \mu(W'_q) \cdot c_i X + O\left( q^2 X^{ \frac{39}{40} } \right) - \mu(q)N_{12}^*(W_q \cap V_{\Z}^{deg,(i)},X)  \right) + \sum_{ q > Q} O_\epsilon \left( (p_1 p_2 \cdots p_n)^{2-\epsilon}  \frac{X}{q^{2-\epsilon} } \right)\\
&&=\sum_{q \leq Q} \left( \mu(W'_q) \cdot c_i X - \mu(q)N_{12}^*(W_q \cap V_{\Z}^{deg,(i)},X)  \right)
+ O_\epsilon \left( (p_1 p_2 \cdots p_n)^{2-\epsilon} X/Q^{1-\epsilon} + X^{\frac{39}{40}}Q^{3+\epsilon} \right)\\
&&= \sum_q c_i \mu(W'_q)X + (p_1 p_2 \cdots p_n)^{2-\epsilon} O_\epsilon \left( X/Q^{1-\epsilon}+X^{\frac{39}{40}}Q^{3+\epsilon} + X^{\frac{199}{200}Q^{1+\epsilon}}\right)\\
&&= c_i \prod_p (1- \mu(W'_q)) X + (p_1 p_2 \cdots p_n)^{2-\epsilon} O_\epsilon \left( X/Q^{1-\epsilon}+X^{\frac{39}{40}}Q^{3+\epsilon} + X^{\frac{199}{200}Q^{1+\epsilon}}\right)\\
&&= c_i \prod_p (\mu(U'_q)) X + (p_1 p_2 \cdots p_n)^{2-\epsilon}O_\epsilon \left( X/Q^{1-\epsilon}+X^{\frac{39}{40}}Q^{3+\epsilon} + X^{\frac{199}{200}Q^{1+\epsilon}}\right).
\end{eqnarray*}

Putting $Q=X^{\frac{1}{400}}$, we have
$$(p_1 p_2 \cdots p_n)^{2-\epsilon}O_\epsilon \left( X/Q^{1-\epsilon}+X^{\frac{39}{40}}Q^{3+\epsilon} + X^{\frac{199}{200}Q^{1+\epsilon}}\right) \ll_\epsilon (p_1 p_2 \cdots p_n)^{2-\epsilon} X^{\frac{399}{400}+\epsilon}.
$$
 Note that
 \begin{eqnarray*}
 c_i \prod_p (\mu(U'_q)) &=& \prod_{k=1}^n \frac{\mu(U_p(\sigma_{p_k}))}{\mu(U_{p_k})} c_i \prod_p \mu(U_p)=  \prod_{k=1}^n \frac{\mu(U_p(\sigma_{p_k}))}{\mu(U_{p_k})} d_i \prod_p \left( 1+ p^{-2} -p^{-4} -p^{-5} \right).
\end{eqnarray*}

From Lemma 20 in \cite{B08}, we can see that, for $f(p)=p^{-1}+2p^{-2}+2p^{-3}+p^{-4}$,
$$
\frac{\mu(U_p(\sigma))}{\mu(U_p)} = \frac{1/120}{1+f(p)},\:\frac{1/12}{1+f(p)},\:\frac{1/8}{1+f(p)},\:\frac{1/6}{1+f(p)},\: \frac{1/6}{1+f(p)},\:\frac{1/4}{1+f(p)},\:\mbox{and }\frac{1/5}{1+f(p)}
$$
for $\sigma=(11111),(1112),(122),(113),(23),(14),(5),$ respectively, and
\begin{eqnarray*}
\frac{\mu(U_p(\sigma))}{\mu(U_p)}& =& \frac{1/6 \cdot 1/p}{1+f(p)},\:\frac{1/2 \cdot 1/p}{1+f(p)},\:\frac{1/3 \cdot 1/p}{1+f(p)},\:\frac{1/2 \cdot 1/p^2}{1+f(p)},\: \frac{1/2 \cdot 1/p^2}{1+f(p)},\:\frac{1/2 \cdot 1/p^2}{1+f(p)},\frac{1/2 \cdot 1/p^2}{1+f(p)},\\
 & & \frac{1/p^3}{1+f(p)},\:\frac{1/p^3}{1+f(p)},\:\mbox{ and } \frac{1/p^4}{1+f(p)}
\end{eqnarray*}
for $\sigma=(1^2111),(1^212),(1^23),(1^21^21),(2^2 1),(1^3 11),(1^3 2),(1^3 1^2), (1^4 1), (1^5),$ respectively. Hence we have proved the theorem. 

\section{Counting $S_4$ quartic fields with local conditions}\label{S_4}

In \cite{BBP}, Belabas, Bhargava and Pomerance obtained a power saving error term for counting $S_4$ quartic fields. For $i=0,1$,
let $N_4^{(i)}(X)$ be the number of isomorphism classes of $S_4$-quartic fields of signature $(4-2i,i)$ with $|d_K|<X$. Then
\begin{eqnarray*}
N_4^{(i)}(X)= D_i X + O(X^{23/24+\epsilon}),
\end{eqnarray*}
where $D_i=d_i \prod_{p} ( 1+p^{-2} -p^{-3} - p^{-4})$, and
$d_0=\frac{1}{48},d_1=\frac{1}{8},$ and $d_0=\frac{1}{16}$.

Let $V_{\Z}$ be the space of pairs of $(A,B)$ of integral ternary quadratic forms. The group $G_\Z=GL_2(\Z) \times SL_3(\Z)$ acts on $V_{\Z}$ in the following way.
For $g_2 \in GL_2(\Z)$ and $g_3 \in SL_3(\Z)$,
\begin{eqnarray*}
(g_2,g_3)\cdot (A,B)^t=g_2( g_3Ag_3^t, g_3Bg_3^t)^t.
\end{eqnarray*}
Here $g_2\cdot (A,B)$ means that $(a_1 (A,B)^t, a_2 (A,B)^t)$, where $a_i$ is the $i$th row of $g_2$.

There is a canonical bijection between the set of $G_\Z$-equivalence classes of elements $(A,B) \in V_{\Z}$, and the set of isomorphism classes of pairs of $(Q,R)$, where $Q$ is a quartic ring and $R'$ is a cubic resolvent ring of $Q$. (See \cite{B04}.)

A prime $p$ has 11 possible splitting type in an $S_4$ quartic field $K$. They are $(1111),$ $(112),$ $(13),$ $(4),$ $(22),$ $(1^211),$ $(1^21^2),$ $(1^22),$ $(2^2),$ $(1^31),$ and $(1^4)$. As we did in the previous section, we can define $U_p(\sigma)$ and $U_p$, resp. and
$\mu(U_p(\sigma))$'s are computed in Lemma 23, \cite{B04}.

Using their result, Yang \cite{yang} was able to count $S_4$ quartic fields with one local condition with a power saving error term. He showed that
\begin{eqnarray*}
N_4^{(i)}(X, \sigma_p)= \frac{\mu(U_p(\sigma_p))}{\mu(U_p)}D_iX + O\left(p^2X^{\frac{143}{144}+\epsilon}\right).
\end{eqnarray*}

By the same argument in the previous section, we can extend Yang's result to the case of finitely many local conditions:
\begin{eqnarray*}
N_4^{(i)}(X, \{\sigma_{p_k} \}_{k=1,2,\cdots,n} ) = \left( \prod_{k=1}\frac{\mu(U_p(\sigma_{p_k}))}{\mu(U_{p_k})} \right) D_i X +
O\left( (p_1\cdots p_k)^2 X^{\frac{143}{144}+\epsilon}\right).
\end{eqnarray*}

From Lemma 23 in \cite{B04}, we can see that, for $f(p)=p^{-1}+2p^{-2}+p^{-3}$,
$$
\frac{\mu(U_p(\sigma))}{\mu(U_p)} = \frac{1/24}{1+f(p)},\:\frac{1/4}{1+f(p)},\:\frac{1/3}{1+f(p)},\:\frac{1/8}{1+f(p)},\: \mbox{and }\frac{1/4}{1+f(p)}
$$
for $\sigma=(1111),(112),(13),(22),(4),$ respectively, and
$$
\frac{\mu(U_p(\sigma))}{\mu(U_p)} = \frac{1/2 \cdot 1/p}{1+f(p)},\:\frac{1/2 \cdot 1/p}{1+f(p)},\:\frac{1/2 \cdot 1/p^2}{1+f(p)},\:\frac{1/2 \cdot 1/p^2}{1+f(p)},\:\frac{1/p^2}{1+f(p)},\: \mbox{and }\frac{1/p^3}{1+f(p)}
$$
for $\sigma=(1^211),(1^22),(1^21^2),(2^2),(1^31), (1^4),$ respectively.  Hence we have proved

\begin{theorem} Let $N_4^{(i)}(X,\mathcal S)$ be the number of $S_4$ quartic fields of signature $(4-2i,i)$ with $|d_K| < X$, and with the local condition $\mathcal S$. Then
\begin{eqnarray*}
N_4^{(i)}(X,\mathcal S)= |\mathcal S| D_i X+ O_{\epsilon}\left(\left(\prod_{p\in \mathcal S} p \right)^2 X^{\frac{143}{144}+\epsilon }\right).
\end{eqnarray*}
\end{theorem}


\begin{thebibliography}{99}
\bibitem{B} P. Billingsley, Probability and Measure, 3rd ed. 1995, John Wiley \& Sons.
\bibitem{BGHT} T. Barnet-Lamb, D. Geraghty, M. Harris, and R. Taylor,
{\em A family of Calabi-Yau varieties and potential automorphy II}, Publ. Res. Inst. Math. Sci. {\bf 47} (2011), 29�-98.
\bibitem{BBP} K. Belabas, M. Bhargava, and C. Pomerance, {\em Error estimates for the Davenport-Heilbronn theorems}, Duke Math. J. \textbf{153} (2010) no.1, 173--210.
\bibitem{B04} M. Bhargava, {\em Higher composition laws III: The parametrization of quartic rings}, Ann. of Math. \textbf{159} (2004), 1329--1360.
\bibitem{B08} \bysame, {\em Higher composition laws IV: The parametrization of quintic rings}, Ann. of Math. \textbf{167} (2008),  53--94.
\bibitem{B10} \bysame, {\em The density of discriminants of quintic rings and fields}, Ann. of Math. (2), \textbf{172} (2010), no. 3, 1559--1591.
\bibitem{CK} P.J. Cho and H. Kim, {\em $n$-level density of Artin $L$-functions}, preprint.
\bibitem{CK1} \bysame, {\em Low lying zeros of Artin $L$-functions}, to appear in Math. Z.
\bibitem{CDF} B. Conrey, W. Duke, and D. Farmer, {\em The distribution of the eigenvalues of Hecke operators},
Acta Arith. {\bf 78} (1997), no. 4, 405--409.
\bibitem{IK} H. Iwaniec and E. Kowalski, Analytic Number Theory, American Mathematical Society Colloquium Publications, 53. American Mathematical Society, Providence, RI, 2004.
\bibitem{JL} G. James and M. Liebeck, {Representations and Characters of Groups}, Cambridge University Press, 1993.
\bibitem{Kl} N. Klingen, Arithmetical Similarities, Oxford Mathematical Monographs, Oxford Science Publications, The
Clarendon Press, Oxford University Press, New York, 1998.
\bibitem{N} H. Nagoshi, {\em Distribution of Hecke eigenvalues}, Proc. of AMS, {\bf 134} (2006), 3097--3106.
\bibitem{Se} J.P. Serre, {\em R\'epartition asymptotique des valeurs propres de l'op\'erateur de Hecke $T_p$}, J. Amer. Math. Soc. {\bf 10} (1997), 75--102.
\bibitem{Se1} \bysame, {\em Quelques applications du th\'eor\`eme de densit\'e de Chebotarev}, Inst. Hautes \'Etudes Sci. Publ. Math. {\bf 54} (1981), 323--401.
\bibitem{ST} A. Shankar and J. Tsimerman, {\em Counting $S_5$ fields with a power saving error term}, preprint.
\bibitem{TT} T. Taniguchi and F. Thorne, {\em Secondary terms in counting functions for cubic fields}, Duke Math. J. {\bf 162} (2013), 2451--2508.
\bibitem{yang} A. Yang, Distribution Problems associated to Zeta Functions and Invariant Theory, Ph.D. Thesis at Princeton University, 2009.

\end{thebibliography}
\end{document}